\documentclass[12pt]{article}
\usepackage{amssymb, amsfonts, amsmath, amsthm}

\def\barr{\begin{array}}
\def\earr{\end{array}}

\def\0{\leqno}
\def\ov{\overline}

\def\Z{{\rlap{$\kern2pt{\rm Z}$}{\rm Z}\,}}

\def\dd{\displaystyle}

\def\dd{\displaystyle}
\def\1{^{-1}}
\def\barr{\begin{array}}
\def\earr{\end{array}}

\def\<{\left<}
\def\>{\right>}

\title{\bf A generalization of Menon's identity}
\author{Marius T\u arn\u auceanu}
\date{May 5, 2012}

\begin{document}

\maketitle

\begin{abstract}
In this note we give a generalization of the well-known Menon's
identity. This is based on applying the Burnside's lemma to a
certain group action.
\end{abstract}



\section{Introduction}

One of the most interesting arithmetical identity is due to P.K.
Menon \cite{7}.

\bigskip\noindent{\bf Menon's Identity.} {\it For every positive integer
$n$ we have
$$\dd\sum_{a\in U(\mathbb{Z}_n)}{\rm gcd}\left(n,a{-}1\right)=\varphi(n)\,\tau(n)\,,$$where
$U(\mathbb{Z}_n)=\{a\in \mathbb{Z}_n \mid {\rm
gcd}\left(n,a\right)=1\}$, $\varphi$ is the Euler's totient
function and $\tau(n)$ is the number of divisors of $n$.}
\bigskip

This identity has many generalizations derived by several authors
(see, for example, \cite{1}-\cite{6}, \cite{10,11} and \cite{13}-\cite{20}). An
usual technique to prove results of this type is based on the
so-called Burnside's lemma (see \cite{8}) concerning group
actions.

\bigskip\noindent{\bf Burnside's Lemma.} {\it Let $G$ be a finite group acting on a finite set $X$ and let $X^g=\{x\in X\mid g\cdot
x=x\}$, for all $g\in G$. Then the number of distinct orbits is
$$N=\frac{1}{\mid G\mid}\dd\sum_{g\in G}\mid X^g\mid.$$}

The starting point for our discussion is given by the open problem
in the end of Section 2 of \cite{16} that suggests to apply the
Burnside's lemma to the natural action of the group $G$ of upper
triangular matrices contained in $GL_r(\mathbb{Z}_n)$, that is
$$G=\left\{\hspace{-1mm}\left(\hspace{-1mm}\begin{array}{llcl}
 a_{11}&a_{12}&\cdots&a_{1r}\\
 0&a_{22}&\cdots&a_{2r}\\
 \vdots&\vdots&&\vdots\\
 0&0&\cdots&a_{rr}\end{array}\hspace{-2mm}\right)\mid a_{ii}\in U(\mathbb{Z}_n) \hspace{1mm}\forall\hspace{1mm} i{=}\ov{1,r},\, a_{ij}\in \mathbb{Z}_n \hspace{1mm}\forall\hspace{1mm} 1{\leq} i{<}j{\leq}
 r\right\}$$(notice that $\mid G\mid\hspace{1mm}=n^{\frac{r(r-1)}{2}}\varphi(n)^r$), on the set
 $$X=\mathbb{Z}_n^r=\left\{\hspace{-1mm}\left(\hspace{-1mm}\begin{array}{llcl}
 x_1\\
 x_2\\
 \vdots\\
 x_r\end{array}\hspace{-2mm}\right)\mid x_i\in\mathbb{Z}_n \hspace{1mm}\forall\hspace{1mm}
 i{=}\ov{1,r}\right\}.$$Following this idea we obtained an
interesting generalization of the Menon's identity. Denote
$\tau_1(n)=\tau(n)$ and $\tau_i(n)=\dd\sum_{d\mid
n}\tau_{i-1}(d)$, for all $i\geq 2$. We remark that an alternative way to define 
these functions is $\tau_i=\tau_{i-1}*e$, where $*$ denotes the usual Dirichlet convolution and $e(n)=1$ for all $n\in\mathbb{N}$. 
In other words, $\tau_i=e^{(i+1)}$, where $e^{(i+1)}$ is the $(i+1)$th power under the Dirichlet convolution (see e.g. \cite{9,12}). Our main result is:

\bigskip\noindent{\bf Theorem.} {\it For every positive integers
$n$ and $r$ we have
$$\dd\sum_{{a_{ii}\in U(\mathbb{Z}_n),\,i=\ov{1,r}}\atop{a_{ij}\in\mathbb{Z}_n},\,1\leq i<j\leq r}\hspace{1mm}\prod_{k=1}^r\hspace{1mm} d_k=n^{\frac{r(r-1)}{2}}\varphi(n)^r\,\tau_r(n)\,,\0(*)$$where
$$\hspace{5mm}d_k={\rm gcd}\left(n,\frac{na_{1k}}{{\rm gcd}\left(n,a_{11}{-}1,a_{12},...,a_{1k-1}\right)}\,,\frac{na_{2k}}{{\rm gcd}\left(n,a_{22}{-}1,a_{23},...,a_{2k-1}\right)}\,,...\,,\right.$$
$$\hspace{-20mm}\left.\frac{na_{k-1k}}{{\rm gcd}\left(n,a_{k-1k-1}{-}1\right)}\,,a_{kk}{-}1\right) \forall\hspace{1mm}k=\ov{1,r}.$$}

\section{Proof of the main theorem}

We will proceed by induction on $r$. Obviously, for $r=1$ the
equality $(*)$ is the Menon's identity.
\bigskip

In the following we will focus on the case $r=2$ (this is not
necessary, but very suggestive for the general implication step).
We have to prove that
$$\dd\sum_{{a_{11},\,a_{22}\in U(\mathbb{Z}_n)}\atop{a_{12}\in\mathbb{Z}_n}}\hspace{-5mm}{\rm gcd}\left(n,a_{11}{-}1\right)\hspace{1mm}{\rm gcd}\left(n,\frac{na_{12}}{{\rm
gcd}\left(n,a_{11}{-}1\right)}\,,a_{22}{-}1\right){=}\,n\varphi(n)^2\,\tau_2(n)\,.\0(1)$$Clearly,
two elements $x=\left(\hspace{-1mm}\begin{array}{llcl}
 x_1\\
 x_2\end{array}\hspace{-2mm}\right)$ and $y=\left(\hspace{-1mm}\begin{array}{llcl}
 y_1\\
 y_2\end{array}\hspace{-2mm}\right)$ of $X$ are contained in the same orbit
if and only if there is $g=\left(\hspace{-1mm}\begin{array}{llcl}
 a_{11}&\hspace{-2mm}a_{12}\\
 0&\hspace{-2mm}a_{22}\end{array}\hspace{-2mm}\right)\in G$ such
that $y=g\cdot x$, i.e.
$$\left\{\begin{array}{llcl}
a_{11}x_1+a_{12}x_2=y_1\\
a_{22}x_2=y_2.\end{array}\right.\0(2)$$We observe that (2) is
equivalent to
$$\left\{\begin{array}{llcl}
\langle x_2\rangle=\langle y_2\rangle\hspace{2mm}(=H)\\
\langle x_1H\rangle=\langle y_1H\rangle\hspace{2mm}\mbox{in }
\mathbb{Z}_n/H,\end{array}\right.\0(2')$$that means
$$\left\{\begin{array}{llcl}
o(x_2)=o(y_2)=\delta\in L_n\\
o_H(x_1)=o_H(y_1)=\delta'\in
L_{\frac{n}{\delta}},\end{array}\right.\0(2'')$$where for a
positive integer $m$ we denote by $L_m$ the lattice of divisors of
$m$. In this way, one obtains
$$N=\hspace{1mm}\mid\{(\delta,\delta')\mid \delta\in L_n, \delta'\in L_{\frac{n}{\delta}}\}\hspace{-1mm}\mid\hspace{1mm}=\dd\sum_{\delta\mid n}\tau(\frac{n}{\delta})=\dd\sum_{\delta\mid n}\tau(\delta)=\tau_2(n).\0(3)$$

\noindent{\bf Remark.} For $r=2$ an explicit formula of $N$ can be
given, namely if $n=p_1^{\alpha_1}p_2^{\alpha_2}\cdots
p_s^{\alpha_s}$ is the decomposition of $n$ as a product of prime
factors, then
$$N=\frac{1}{2^s}\prod_{i=1}^s (\alpha_i+1)(\alpha_i+2).$$
\bigskip

Next we observe that for a fixed
$g=\left(\hspace{-1mm}\begin{array}{llcl}
 a_{11}&\hspace{-2mm}a_{12}\\
 0&\hspace{-2mm}a_{22}\end{array}\hspace{-2mm}\right)\in G$ we have $x=\left(\hspace{-1mm}\begin{array}{llcl}
 x_1\\
 x_2\end{array}\hspace{-2mm}\right)\in X^g$ if and only
 if $g\cdot x=x$, i.e.
$$\left\{\begin{array}{llcl}
(a_{11}-1)x_1+a_{12}x_2=0\\
(a_{22}-1)x_2=0.\end{array}\right.\0(4)$$By multiplying the first
equation with $\frac{n}{{\rm gcd}(n,\,a_{11}-1)}$\,, it follows
that (4) is e\-qui\-va\-lent to
$$\left\{\begin{array}{llcl}
\frac{na_{12}}{{\rm gcd}(n,\,a_{11}-1)}\,\,x_2=0\\
(a_{22}-1)x_2=0\end{array}\right.\0(4')$$and consequently to
$$x_2{\in}\langle\frac{n}{{\rm gcd}(n,\,\frac{na_{12}}{{\rm gcd}(n,\,a_{11}-1)})}\rangle\cap\langle\frac{n}{{\rm gcd}(n,\,a_{22}-1)}\rangle{=}\langle\frac{n}{{\rm gcd}(n,\frac{na_{12}}{{\rm
gcd}(n,\,a_{11}-1)}\,,a_{22}-1)}\rangle.$$So, $x_2$ can be chosen
in ${\rm gcd}(n,\frac{na_{12}}{{\rm
gcd}(n,\,a_{11}-1)}\,,a_{22}-1)$ ways. Moreover, we easily infer
that for each such choice $x_1$ can be chosen in ${\rm
gcd}(n,\,a_{11}-1)$ ways. Hence
$$\mid X^g\mid\hspace{1mm}={\rm gcd}(n,\,a_{11}-1){\rm gcd}\left(n,\frac{na_{12}}{{\rm
gcd}\left(n,\,a_{11}-1\right)}\,,a_{22}-1\right),$$which together
with (3) lead to (1), as desired.
\bigskip

Finally, we will prove the general implication step. Assume that
$(*)$ holds for $r-1$. Two elements
$x=\left(\hspace{-1mm}\begin{array}{llcl}
 x_1\\
 x_2\\
 \vdots\\
 x_r\end{array}\hspace{-2mm}\right)$ and $y=\left(\hspace{-1mm}\begin{array}{llcl}
 y_1\\
 y_2\\
 \vdots\\
 y_r\end{array}\hspace{-2mm}\right)$ of $X$ belong to the same orbit
if and only if $y=g\cdot x$ for some
$g=\left(\hspace{-1mm}\begin{array}{llcl}
 a_{11}&a_{12}&\cdots&a_{1r}\\
 0&a_{22}&\cdots&a_{2r}\\
 \vdots&\vdots&&\vdots\\
 0&0&\cdots&a_{rr}\end{array}\hspace{-2mm}\right)\in G$, i.e.
$$\left\{\begin{array}{llcl}
a_{11}x_1+a_{12}x_2+\cdots +a_{1r}x_r=y_1\\
a_{22}x_2+a_{23}x_3+\cdots +a_{2r}x_r=y_2\\
\vdots\\
a_{rr}x_r=y_r.\end{array}\right.\0(5)$$The equalities (5) are
equivalent to
$$\left\{\begin{array}{llcl}
\langle x_r\rangle=\langle y_r\rangle\\
\langle x_{r-1}H_1\rangle=\langle
y_{r-1}H_1\rangle\hspace{2mm}\mbox{in } \mathbb{Z}_n/H_1\\
\vdots\\
\langle x_1H_{r-1}\rangle=\langle
y_1H_{r-1}\rangle\hspace{2mm}\mbox{in } \mathbb{Z}_n/H_{r-1},
\end{array}\right.\0(5')$$where
\begin{itemize}
\item[] \hspace{30mm}$H_1\hspace{-0,5mm}=\hspace{-0,5mm}\langle x_r\rangle\hspace{-0,5mm}=\hspace{-0,5mm}\langle y_r\rangle$,
\item[] \hspace{30mm}$H_2\hspace{-0,5mm}=\hspace{-0,5mm}\langle x_{r-1},x_r\rangle\hspace{-0,5mm}=\hspace{-0,5mm}\langle y_{r-1},y_r\rangle$,
\item[] \hspace{30mm}\vdots
\item[] \hspace{30mm}$H_{r-1}\hspace{-0,5mm}=\hspace{-0,5mm}\langle x_2,x_3,...,x_r\rangle\hspace{-0,5mm}=\langle y_2,y_3,...,y_r\rangle$,
\end{itemize}
\noindent which means
$$\left\{\begin{array}{llcl}
o(x_r)=o(y_r)=\delta_1\in L_n\\
o_{H_1}(x_{r-1})=o_{H_1}(y_{r-1})=\delta_2\in
L_{\frac{n}{\delta_1}}\\
\vdots\\
o_{H_{r-1}}(x_1)=o_{H_{r-1}}(y_1)=\delta_r\in
L_{\frac{n}{\delta_1\delta_2\cdots\delta_{r-1}}}\,.
\end{array}\right.\0(5'')$$It is now easy to see that
$$N=\hspace{1mm}\mid\{(\delta_1,\delta_2,...,\delta_{r})\mid \delta_1\in L_n, \delta_2\in L_{\frac{n}{\delta_1}},...,\delta_r\in
L_{\frac{n}{\delta_1\delta_2\cdots\delta_{r-1}}}\}\hspace{-1mm}\mid\hspace{1mm}=\0(6)$$
$$\hspace{2mm}=\dd\sum_{\delta_1\mid n}\mid\{(\delta_2,...,\delta_{r})\mid \delta_2\in L_{\frac{n}{\delta_1}},...,\delta_r\in
L_{\frac{n}{\delta_1\delta_2\cdots\delta_{r-1}}}\}\hspace{-1mm}\mid\hspace{1mm}=\cdots=$$
$$\hspace{-6mm}=\dd\sum_{\delta_1\mid n}\hspace{3mm}\dd\sum_{\delta_2\mid
\frac{n}{\delta_1}}\hspace{3mm}\cdots\hspace{-5mm}\dd\sum_{\hspace{5mm}\delta_{r-1}\mid\frac{n}{\delta_1\delta_2\cdots\delta_{r-2}}}\tau(\frac{n}{\delta_1\delta_2\cdots\delta_{r-1}})=\tau_r(n).$$On
the other hand, given
$g\hspace{-1mm}=\hspace{-1mm}\left(\hspace{-1mm}\begin{array}{llcl}
 a_{11}&a_{12}&\cdots&a_{1r}\\
 0&a_{22}&\cdots&a_{2r}\\
 \vdots&\vdots&&\vdots\\
 0&0&\cdots&a_{rr}\end{array}\hspace{-2mm}\right)\hspace{-1mm}\in G$, we have $x\hspace{-1mm}=\hspace{-1mm}\left(\hspace{-1mm}\begin{array}{llcl}
 x_1\\
 x_2\\
 \vdots\\
 x_r\end{array}\hspace{-2mm}\right)\hspace{-1mm}\in X^g$ if and
only if $$\left\{\begin{array}{llcl}
(a_{11}-1)x_1+a_{12}x_2+\cdots+a_{1r}x_r=0\\
(a_{22}-1)x_2+a_{23}x_2+\cdots+a_{2r}x_r=0\\
\vdots\\
(a_{rr}-1)x_r=0.\end{array}\right.\0(7)$$By multiplying the first
equation with $\frac{n}{{\rm
gcd}(n,\,a_{11}-1,\,a_{12},...,\,a_{1r-1})}$\,, the second one
with $\frac{n}{{\rm
gcd}(n,\,a_{22}-1,\,a_{23},...,\,a_{2r-1})}$\,,..., and the last
but one with $\frac{n}{{\rm gcd}(n,\,a_{r-1r-1}-1)}$\,, (7)
becomes a system in $x_r$ that has $$d_r={\rm
gcd}\left(n,\frac{na_{1r}}{{\rm
gcd}\left(n,a_{11}{-}1,a_{12},...,a_{1r-1}\right)}\,,...,\,\frac{na_{r-1r}}{{\rm
gcd}\left(n,a_{r-1r-1}{-}1\right)}\,,a_{rr}{-}1\right)$$solutions,
namely $x_r\in \langle\frac{n}{d_r}\rangle$. Put
$x_r=\gamma\frac{n}{d_r}$ with $\gamma\in\{0,1,...,d_r-1\}$. Then
(7) can be rewritten as
$$\left\{\begin{array}{llcl}
(a_{11}-1)x_1+a_{12}x_2+\cdots+a_{1r-1}x_{r-1}=-\gamma\frac{n}{d_r}a_{1r}\\
(a_{22}-1)x_2+a_{23}x_3+\cdots+a_{2r-1}x_{r-1}=-\gamma\frac{n}{d_r}a_{2r}\\
\vdots\\
(a_{r-1r-1}-1)x_{r-1}=-\gamma\frac{n}{d_r}a_{r-1r}.\end{array}\right.\0(7')$$If
$(x_1^0,x_2^0,...,x_{r-1}^0)$ is a particular solution of (7'),
then one obtains a homogeneous system
$$\left\{\begin{array}{llcl}
(a_{11}-1)(x_1-x_1^0)+a_{12}(x_2-x_2^0)+\cdots+a_{1r-1}(x_{r-1}-x_{r-1}^0)=0\\
(a_{22}-1)(x_2-x_2^0)+a_{23}(x_3-x_3^0)+\cdots+a_{2r-1}(x_{r-1}-x_{r-1}^0)=0\\
\vdots\\
(a_{r-1r-1}-1)(x_{r-1}-x_{r-1}^0)=0\end{array}\right.\0(7'')$$with
$\prod_{k=1}^{r-1}\hspace{1mm} d_k$ solutions by the inductive
hypothesis. We infer that
$$\mid X^g\mid\hspace{1mm}=\prod_{k=1}^r\hspace{1mm} d_k,$$which
together with (6) lead to the equality $(*)$. This completes the
proof. \hfill\rule{1,5mm}{1,5mm}

\bigskip\noindent{\bf Acknowledgements.} The author is grateful to the reviewer for
its remarks which improve the previous version of the paper.

\vspace*{5ex}\small

\hfill
\begin{minipage}[t]{5cm}
Marius T\u arn\u auceanu \\
Faculty of  Mathematics \\
``Al.I. Cuza'' University \\
Ia\c si, Romania \\
e-mail: {\tt tarnauc@uaic.ro}
\end{minipage}

\end{document}